\documentclass[a4paper,12pt]{article}

\pdfoutput=1

\usepackage[paper=letterpaper,margin=1in]{geometry}
\usepackage{amsrefs,amsmath,amssymb,amsfonts,epsfig,cite,setspace,bigstrut,longtable,array,breqn,color}
\usepackage{caption,setspace}
\usepackage{amsthm}
\usepackage{pdfpages,lipsum}
\usepackage{todonotes}
\usepackage{cancel}
\usepackage[export]{adjustbox}

\usepackage[OT2,T1]{fontenc}
\DeclareSymbolFont{cyrletters}{OT2}{wncyr}{m}{n}
\DeclareMathSymbol{\Sha}{\mathalpha}{cyrletters}{"58}

\usepackage{mathtools}
\usepackage{hyperref}

\usepackage{listings}
\lstdefinelanguage{Sage}[]{Python}
{morekeywords={False,sage,True},sensitive=true}
\definecolor{dblackcolor}{rgb}{0.0,0.0,0.0}
\definecolor{dbluecolor}{rgb}{0.01,0.02,0.7}
\definecolor{dgreencolor}{rgb}{0.2,0.4,0.0}
\definecolor{dgraycolor}{rgb}{0.30,0.3,0.30}

\newenvironment{red}{\relax\color{red}}{\relax}
\newenvironment{blue}{\relax\color{blue}}{\hspace*{.5ex}\relax}
\newcommand{\ber}{\begin{red}}
\newcommand{\er}{\end{red}}
\newcommand{\beb}{\begin{blue}}
\newcommand{\eb}{\end{blue}}

\begin{document}


\vskip 0.25in

\newcommand{\nn}{\nonumber}
\newcommand{\tr}{\mathop{\rm Tr}}
\newcommand{\comment}[1]{}
\newcommand{\cM}{{\cal M}}
\newcommand{\cW}{{\cal W}}
\newcommand{\cN}{{\cal N}}
\newcommand{\cH}{{\cal H}}
\newcommand{\cK}{{\cal K}}
\newcommand{\cZ}{{\cal Z}}
\newcommand{\cO}{{\cal O}}
\newcommand{\cA}{{\cal A}}
\newcommand{\cB}{{\cal B}}
\newcommand{\cC}{{\cal C}}
\newcommand{\cD}{{\cal D}}
\newcommand{\cT}{{\cal T}}
\newcommand{\cV}{{\cal V}}
\newcommand{\cE}{{\cal E}}
\newcommand{\cF}{{\cal F}}
\newcommand{\cX}{{\cal X}}
\newcommand{\IA}{\mathbb{A}}
\newcommand{\IP}{\mathbb{P}}
\newcommand{\IQ}{\mathbb{Q}}
\newcommand{\IH}{\mathbb{H}}
\newcommand{\IR}{\mathbb{R}}
\newcommand{\IC}{\mathbb{C}}
\newcommand{\IF}{\mathbb{F}}
\newcommand{\IV}{\mathbb{V}}
\newcommand{\II}{\mathbb{I}}
\newcommand{\IZ}{\mathbb{Z}}
\newcommand{\re}{{\rm~Re}}
\newcommand{\im}{{\rm~Im}}

\let\oldthebibliography=\thebibliography
\let\endoldthebibliography=\endthebibliography
\renewenvironment{thebibliography}[1]{%
\begin{oldthebibliography}{#1}%
\setlength{\parskip}{0ex}%
\setlength{\itemsep}{0ex}%
}%
{%
\end{oldthebibliography}%
}

\newtheorem{theorem}{\bf THEOREM}
\def\thetheorem{\thesection.\arabic{theorem}}
\newtheorem{proposition}{\bf PROPOSITION}
\def\thetheorem{\thesection.\arabic{proposition}}
\newtheorem{observation}{\bf OBSERVATION}
\def\thetheorem{\thesection.\arabic{observation}}
\newtheorem{conjecture}{\bf CONJECTURE} 
\def\thetheorem{\thesection.\arabic{CONJECTURE}}

\theoremstyle{definition}
\newtheorem{definition}{\bf DEFINITION} 
\def\thetheorem{\thesection.\arabic{DEFINITION}}
\newtheorem{example}{\bf EXAMPLE} 
\def\thetheorem{\thesection.\arabic{EXAMPLE}}
\newtheorem{remark}{\bf REMARK} 
\def\thetheorem{\thesection.\arabic{REMARK}}

\def\theequation{\thesection.\arabic{equation}}
\newcommand{\setall}{\setcounter{equation}{0}
        \setcounter{theorem}{0}}
\newcommand{\setequation}{\setcounter{equation}{0}}
\renewcommand{\thefootnote}{\fnsymbol{footnote}}

\newcommand{\seteq}{\mathbin{:=}}
\newcommand{\GL}{\operatorname{GL}}
\newcommand{\Sp}{\operatorname{Sp}}
\newcommand{\USp}{\operatorname{USp}}
\newcommand{\GSp}{\operatorname{GSp}}
\newcommand{\SU}{\operatorname{SU}}
\newcommand{\SO}{\operatorname{SO}}
\newcommand{\End}{\operatorname{End}}

\begin{titlepage}

~\\
\vskip 1cm

\begin{center}
{\Large \bf Murmurations of Elliptic Curves}
\end{center}
\medskip

\renewcommand{\arraystretch}{0.5} 

\vspace{.4cm}
\centerline{
{\large Yang-Hui He, Kyu-Hwan Lee, Thomas Oliver, Alexey Pozdnyakov}
}
\vspace*{3.0ex}

\vspace{10mm}

\begin{abstract}

We investigate the average value of the Frobenius trace at $p$ over elliptic curves in a fixed conductor range with given rank.
Plotting this average as $p$ varies over the primes yields a striking oscillating pattern, the details of which vary with the rank.
Based on this observation, we perform various data-scientific experiments with the goal of classifying elliptic curves according to their ranks.

\end{abstract}
\begin{center}
\includegraphics[scale=0.5]{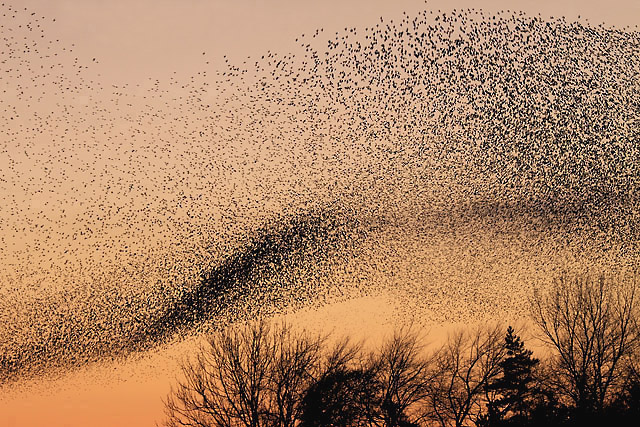}
\end{center}
(Photo by Walter Baxter / {\em A murmuration of starlings at Gretna} / CC BY-SA 2.0)
\end{titlepage}

\begin{spacing}{1}
\tableofcontents
\end{spacing}

\section{Introduction}
Elliptic curves are important objects in number theory, not only in their own right, but also for the crucial role they have played in a range of developments. 
A significant example is Fermat's Last Theorem, which was a consequence of the {\em modularity of elliptic curves} conjectured by Taniyama, Shimura, and Weil, and established by  Andrew Wiles \cite{W95}.  
Despite extensive studies, there are still many open problems in the theory of elliptic curves, the most celebrated example of which is perhaps the Birch and Swinnerton-Dyer (BSD) conjecture.

The BSD conjecture is concerned with the ranks of elliptic curves. 
If $E$ is an elliptic curve over $\mathbb Q$, then, by the Mordell--Weil theorem, we know that the rational points of $E$ form a finitely generated abelian group. 
Though the torsion part of the abelian group is relatively well understood, the rank of the free part is still mysterious. 
In particular, there is no general algorithm to compute the rank of an elliptic curve, and we do not know whether or not it can be arbitrarily large. 

The BSD conjecture relates the rank of an elliptic curve (an aspect of its algebraic structure) to an analytic property of its $L$-function.  
For $\mathrm{Re}(s)\gg0$, the $L$-function of an elliptic curve over $\mathbb{Q}$ may be written as an Euler product \[L(E,s)=\prod_{p\text{ prime}}L_p(E,s)^{-1},\]
where $L_p(E,s)\in\mathbb{Z}[p^{-s}]$ has degree $\leq2$.
If $E$ has good reduction at a prime $p$, then
\[L_p(E,s)=1-a_p(E)p^{-s}+p^{1-2s},\] 
where $a_p(E)=p+1-\#E(\mathbb F_p)$, and $\#E(\mathbb F_p)$ is the number of points of $E$ over $\mathbb F_p$.

Many fundamental invariants of $E$ are connected to analytic properties of $L(E,s)$. 
For example, the conductor of $E$ appears in the functional equation for $L(E,s)$. 
Most importantly, the (weak) BSD conjecture predicts that the rank of $E$ is equal to the order of vanishing for $L(E,s)$ at $s=1$.

In this paper, we study the average value of $a_p(E)$, as $E$ varies over the set of elliptic curves over $\mathbb{Q}$ with fixed rank and conductor in a specified range.
In this introduction we provide a flavour of the resulting images from Section~\ref{sec:av}, with complete details given therein.

Enumerating the primes in ascending order yields the sequence $p_1=2,~p_2=3,~p_3=5, \dots$.
Given an elliptic curve $E$, we will refer to the sequence $(a_{p_i}(E))_{i=1}^{\infty}$ as the $a_p$-coefficients for $E$.
We are interested in the following function of $n \in\mathbb{Z}_{\geq1}$:
\begin{equation}\label{eq.frp}
f_r(n)=\frac{1}{\#\mathcal{E}_r[N_1,N_2]}\sum_{E\in\mathcal{E}_r[N_1,N_2]}a_{p_n}(E),
\end{equation}
where $N_1 < N_2\in\mathbb{Z}_{>0}$, and $\mathcal{E}_r[N_1,N_2]$ is the set of (representatives for the isogeny classes of) elliptic curves over $\mathbb{Q}$ with rank $r$ and conductor in range $[N_1,N_2]$.

\begin{example}\label{ex.firstplots}
We have $p_{1000}=7919$.
Using \cite[Elliptic~curves~over~$\mathbb{Q}$]{lmfdb}, we see that 
\[\# \mathcal E_0[7500,10000]=4328, \ \ \# \mathcal E_1[7500,10000]=5194, \] 
\[  \# \mathcal E_0[5000, 10000]=8536, \ \  \# \mathcal E_2[5000, 10000]=1380.\] 
For $1 \le n \le 1000$ (i.e., for primes $2\leq p\leq7919$), plotting the points $(n,f_r(n))$ for $r \in \{0,1\}$ and $[N_1,N_2]=[7500,10000]$ (resp. $r \in \{0,2\}$ and $[N_1,N_2]=[5000,10000]$), yields the top (resp. bottom) image in Figure \ref{f-012}, in which blue and red (resp. blue and green) dots represent $f_0(n)$ and $f_1(n)$ (resp. $f_0(n)$ and $f_2(n)$), respectively.
\begin{figure}[t!!!]
\begin{center}
     \includegraphics[width=0.9\textwidth]{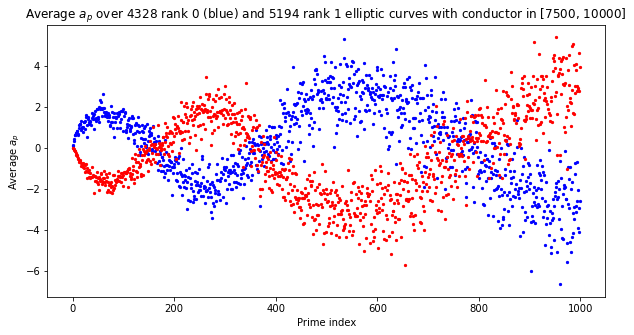}
     \includegraphics[width=0.9\textwidth]{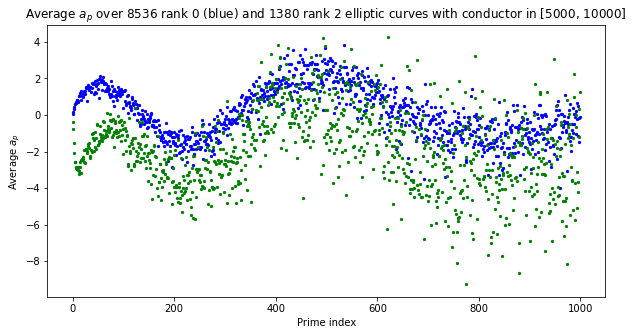}
\end{center}
\caption{\sf (Top) Plots of the functions $f_0(n)$ (blue) and $f_1(n)$ (red) for $1 \le n \le 1000$ and $[N_1,N_2]=[7500,10000]$. (Bottom) Plots of the functions $f_0(n)$ (blue) and $f_2(n)$ (green) for $1 \le n \le 1000$ and $[N_1,N_2]=[5000,10000]$. Further details are given in Example~\ref{ex.firstplots}.}
\label{f-012}
\end{figure}
\end{example}

To the best knowledge of the authors, this oscillating behaviour for the average values of $a_p$ has never been reported in the literature. 
Moreover, varying the rank yields strikingly distinctive patterns, which we believe may be exploited to advance our understanding of the ranks of elliptic curves. 
 
There is nothing special about the intervals $[N_1,N_2]$ used in Example~\ref{ex.firstplots}. 
The only requirement is that they are neither too narrow nor too wide compared to the primes $p$. 
Indeed, as shown in Section~\ref{s:fitting}, similar patterns may be observed for other choices.
Moreover, as observed by Sutherland, there is a certain invariance to the oscillations as the intervals are scaled to include larger conductors \cite{Suth}. 
This so-called {\em scale invariance} is treated systematically in \cite{HLOPS}, which also explores the average values of $a_p$-coefficients attached to modular forms and genus $2$ curves.

 Though the average in equation~\eqref{eq.frp} was first explored from a data scientific perspective -- wherein the primes were simply a list of features indexed by $n$ -- the overwhelming weight of theoretical and experimental evidence suggests an underlying mathematical structure.
From an arithmetic perspective, it is more natural to view the average as a function of the primes $p$, rather than the prime index $n$, though the prime number theorem implies that $p_n$ is asymptotic to $n\log(n)$ and so this distinction yields only a logarithmic difference. 
With this in mind, in Section~\ref{s:fitting}, we consider the functions 
\begin{equation}
g_r(p)=\frac{1}{\#\mathcal{E}_r[N_1,N_2]}\sum_{E\in\mathcal{E}_r[N_1,N_2]}a_{p}(E),
\end{equation}
and investigate the oscillations of $g_r(p)$ for several conductor ranges. 
We then fit the oscillations with curves of the form

\[y=Ax^{\alpha}\sin\left(Bx^{\beta}\right)\] 
by determining numerical values for the parameters $A,B,\alpha,\beta$ which minimize the mean squared error.

This paper is a continuation of the recent work \cite{HLOa}, \cite{HLOb}, \cite{HLOc} by the first three authors, where they applied machine-learning techniques to distinguish arithmetic curves and number fields according to standard invariants such as rank, Sato--Tate group, class number, and Galois group. 
The experimental results clearly show that these number theoretic objects can be classified by machine-learning with high accuracy (often $> 97\%$) once they are presented appropriately, and demonstrate the capacity of machine-learning to predict basic invariants of objects in algebraic number theory.
In particular, in \cite{HLOc}, an elliptic curve $E$ of rank $0$ or $1$ is presented by the $500$-dimensional vector whose $n$th co-ordinate is $a_{p_n}(E)$.
Vectors corresponding to curves of different ranks were distinguished by logistic regression with accuracy $> 97\%$. 

As was pointed out by several experts after \cite{HLOc} was posted, the parity conjecture implies that the distinction between rank $0$ and rank $1$ could be made by observing the sign (i.e. root number) in the functional equation. 
Whilst this is true, it is not clear how one can compute the root number from only finitely many $a_p(E)$. 
Moreover, parity does not distinguish curves between rank $0$ from rank $2$. In this paper, we show that the distinction of rank $0$ and rank $2$ can also be made through logistic regression with high accuracy, and apply PCA to see clustering of elliptic curves of rank $0$, $1$ and $2$. 
One can see that Figure \ref{f-012} already explains why logistic regression is so efficient. 
Namely, for small $p$, the values of $a_p(E)$ behave quite differently (on average), depending on the rank of $E$. 

We conclude this introduction with an overview of what is to come. 
In Section~\ref{s:strategy} we review some basic theory of elliptic curves. 
In Section~\ref{s:datasets}, we describe some data-scientific concepts utilised in the sequel, for example, point clouds, logistic regression, and principal component analysis (PCA).
In Section~\ref{s:elliptic} we describe our experimental results, including curve fitting.
In Section~\ref{sec:heuristic}, we develop a heuristic formula for distinguishing curves of rank $<2$ from those of rank $\geq2$. 
In Section~\ref{s:outlook}, we discuss possible extensions of our work, both theoretical and experimental. 

\subsection*{Acknowledgements}
All the elliptic curves used for this paper are downloaded from LMFDB \cite{lmfdb} and their $a_p$-coefficients are generated by SageMath \cite{sage}. 
We thank Andrew Sutherland for numerous helpful comments on a preliminary draft, 
Minhyong Kim for several insights, 
and Chris Wuthrich for useful discussions, and the anonymous referees for their suggested revisions.
YHH is indebted to STFC UK, for grant ST/J00037X/2, KHL is partially supported by a grant from the Simons Foundation (\#712100), and TO acknowledges support from the EPSRC through research grant EP/S032460/1.

\section{Elliptic curves}\label{s:strategy}

 In this section we review the necessary mathematical background. 
Let $X$ be a smooth, projective, geometrically connected curve of genus $g$ defined over $\mathbb{Q}$. 
We say that a prime number $p$ is a {\em good prime} for $X$ if there exists an integral model for $X$ whose reduction modulo $p$ defines a smooth variety of the same dimension.
In this work, we will focus on the case that $X=E$ is an elliptic curve defined over $\mathbb{Q}$. 
In particular, we have $g=1$.
If $p$ is a good prime for $E$, then we introduce the polynomial:
\begin{equation}\label{eq.EllipticEuler}
L_p(E,T)=1 - a_p(E)T + p T^2 , \ \ (p\text{ good}),
\end{equation}
where
\begin{equation}\label{eq.ap}
a_p(E)=p+1-\#E\left(\mathbb{F}_{p}\right).
\end{equation}
If $p$ is a bad prime for $E$, the analogue of equation~\eqref{eq.EllipticEuler} is $L_p(E,T)=1-a_p(E)T$ with $a_p(E) \in \{-1,0,1\}$, depending on the reduction type.
If $E$ is given in Weierstrass form, then equation~\eqref{eq.ap} is valid for all $p$, good or bad.
More information can be found, e.g., in \cite[Appendix C.16]{Sil1}. The $L$-function of $E$ is defined by 
\[L(E,s)=\prod_{p\text{ prime}}L_p(E,p^{-s})^{-1}.\] 
By the Mordell--Weil Theorem, the set of rational points $E(\mathbb Q)$ of an elliptic curve $E$ defined over $\mathbb{Q}$ forms a finitely generated abelian group and thus decomposes into a product of the torsion part $E(\mathbb Q)_{\mathrm{tor}}$ and the free part:
\[   E(\mathbb{Q}) \cong E(\mathbb Q)_{\mathrm{tor}} \times \mathbb Z^{r_E}. \] 
The rank $r_E$ of the free part is called the {\em rank} of the elliptic curve $E$.

The rank has been a focal point of extensive studies on elliptic curves, but there still remain important open problems. 
In particular, we do not know whether the set of ranks is bounded or not.

In fact, the largest known-rank, established by Noam Elkies \cite{E07} in 2006, is (expected to be) 28, and the boundedness of the average rank of elliptic curves was proved by Bhargava and Shankar \cite{BS1, BS2} in 2015. 
As mentioned in the introduction, the celebrated Birch and Swinnerton-Dyer (BSD) conjecture predicts that the rank of $E$ is equal to the order of vanishing of $L(E,s)$ at $s=1$. 
For rank 0 and 1, this conjecture is known to be true  by the work of Kolyvagin \cite{K89} and the modularity theorem \cite{W95, BCDT}.

As for the Riemann $\zeta$-function, the $L$-function $L(E,s)$ may be completed to a function $\Lambda(E,s) = \pi^{-s} \Gamma(\frac s 2) \Gamma(\frac {s+1} 2) N_E^{s/2} L(E,s)$ which admits analytic continuation to $\mathbb C$ and satisfies a functional equation
\[  \Lambda (E,2- s) = w(E) \Lambda(E,2), \] 
where $N_E \in \mathbb Z$ is the  conductor of $E$ (see \cite[Section~VII.11]{Sil1}) and $w(E)= \pm 1$ is the root number. 
Assuming that the BSD conjecture is true, the {\em parity conjecture} asserts that:
\begin{equation}\label{eq.parity}
    (-1)^{r_E} = w(E) .
\end{equation} 
For elliptic curves of ranks 0 and 1, Equation~\eqref{eq.parity} is a theorem.
Thus, if $r_E\in\{0,1\}$, then we can determine the rank $r_E$ by looking at the root number $w(E)$. 
However, the same argument would not work, for example, if $r_E\in\{0,2\}$. 
\section{Datasets and strategies}\label{s:datasets}

In this section, we explain how to make our datasets of elliptic curves and give an overview of the machine-learning strategies used.

\subsection{Point clouds of elliptic curves}\label{s.pointclouds}
For $i\in\mathbb{Z}_{>0}$, let $p_i$ denote the $i$th prime. 
In particular, $p_1=2$ and $p_{1000}=7919$. 
For an elliptic curve $E$, we introduce the vector:
\begin{equation}\label{eq.Lvector}
v_L(E)=\left(a_{p_1}(E),\dots,a_{p_{1000}}(E)\right)\in\mathbb{Z}^{1000},
\end{equation}
where $a_p(E)$ is defined in \eqref{eq.ap}. 
 
One could implement what follows with any dimension $d$ in place of $1000$.
For example, in \cite{HLOc}, the first three authors implemented similar experiments with $d\in\{100,200,300,500\}$\footnote{ We caution the reader that, in \cite{HLOc}, the letter $N$ is used instead of $d$. In the present text, the letter $N$ is reserved for the conductor of an elliptic curve.}.
Furthermore, in Section \ref{sec:heuristic}, we will develop classifiers in much lower dimensions (namely $d=10$).
Naively, one might expect that using a larger value for $d$ may increase the accuracy of the classifiers, though, as shown in \cite{MMNS}, this matter is somewhat subtle. As for the oscillatory behavior of the functions $f_r(n)$ defined in \eqref{eq.frp}, we need to have $d$ reasonably big and comparable with the conductor ranges. See figures in Section \ref{s:elliptic}. 

In our dataset, an elliptic curve $E$ is represented by the vector $v_L(E) \in \mathbb{R}^{1000}$ and we study the properties of the collection of vectors $\{v_L(E)\}_E\subset\mathbb{R}^{1000}$. 
In the parlance of data science, we investigate $\{ v_L(E) \}$ as a point cloud.
Each point may be further labelled with properties of $E$ such as its rank $r_E$, or conductor $N_E$.
Various investigations of machine-learning on this dataset was performed in \cite{HLOc}, including classification of elliptic curves of rank 0 and 1.

We use Cremona's database \cite{Cre} of elliptic curves, which can also be accessed through \cite{lmfdb}. 
The completeness of the database is discussed in \cite[Completeness~of~elliptic~curve~data~over~$\mathbb{Q}$]{lmfdb}. 

Using \cite{sage}, $a_{p_n}(E)$ can be calculated for each $E$ and $n$. For example, when $E$ is given by label `37a1', it is defined by   $y^2+y=x^3-x$, and we can get its rank and $a_{p_n}$ for $n=15$ as follows:  
\begin{lstlisting}
sage: E=EllipticCurve('37a1'); n=15
sage: r=E.rank()
sage: v=E.aplist(Primes()[n-1])
sage: E, r, v
(Elliptic Curve defined by y^2 + y = x^3 - x over Rational Field,
 1,
 [-2, -3, -2, -1, -5, -2, 0, 0, 2, 6, -4, -1, -9, 2, -9])
\end{lstlisting}
More details can be found in  \cite[Elliptic curves over the rational numbers]{ELS}. In this way, we obtain our datasets consisting of $v_L(E)$ labelled according to rank for various ranges of conductors $N_E$.

We note that the Hasse--Weil $L$-function $L(E,s)$ of an elliptic curve $E$ is an invariant of its isogeny class, and our datasets actually have only one representative curve for each isogeny class.

\subsection{Averaging}
The size of the set $\mathcal{E}_r[N_1,N_2]$ varies with the parameters $r$, $N_1$, and $N_2$.
In this article, we will choose the parameters so that $\mathcal{E}_r[N_1,N_2]$ has order approximately $k\times10^3$ for $1<k<10$.
By averaging (arithmetic mean), we construct a single value $f_r(n)$, representative of the set of values $\{a_{p_n}(E):E\in\mathcal{E}_r[N_1,N_2]\}$.
It seems unlikely that any genuine elliptic curve could have $a_{p_n}$ staying near $f_r(n)$ for all $n$. 
One may consider the geometric mean; however, we do not observe any interesting features in its distribution.
 The standard deviation of $a_{p_n}$ is asymptotically equal to $\sqrt{p_n}$ for any $r$ and will not play any role in our discussion.

\subsection{Logistic regression}\label{s:logreg}

The binary logistic regression classifier is a strategy for supervised machine-learning based upon the logistic sigmoid function
\[\sigma:\mathbb{R}\rightarrow(0,1), \ \ \sigma(x)=\frac{1}{1+e^{-x}}.\]
In multi-class logistic regression, its generalization, called the softmax function, is used.
Further details may be found in \cite[Sections~4.4,~11.3]{Hastie}.

In Section~\ref{sec:ECrank} we present various binary and ternary experiments involving elliptic curves of rank $r_E\in\{0,1,2\}$.
Recall from Section~\ref{s.pointclouds} that each elliptic curve defines a vector $v_L(E)\in\mathbb{R}^{1000}$. 
A binary logistic regression classifier works by finding a single vector $\mathbf w\in\mathbb{R}^{1000}$ and number $b\in\mathbb{R}$ so that 
\[\sigma(v_L(E) \cdot  \mathbf w+b)\] 
is a predictor for the rank $r_E$ of $E$ where $v_L(E)\cdot  \mathbf w\in\mathbb{R}$ denotes the dot product of $v_L(E)$ and $\mathbf w$. 
In each experiment, the vectors $\mathbf w$ and numbers $b$ are calculated by numerical means, and we do not make them explicit. 
The multi-class case is similar with $\sigma$ replaced by the softmax function.

An explicit binary logistic regression experiment is presented in Section~\ref{sec:heuristic}, and involves two sets of elliptic curves with conductors in a specified range: those with rank $r_E<2$ and those with rank $r_E\geq2$. 
This time, we present each elliptic curve by the $10$-dimensional vector $\mathbf a=\mathbf a (E)=(a_{p_1},\dots,a_{p_{10}})=(a_2,\dots,a_{29})\in\mathbb{R}^{10}$ (a projection of $v_L(E)$).
We find an explicit vector $\mathbf w\in\mathbb{R}^{10}$ and $b \in \mathbb R$ such that 
\[\sigma(\mathbf a \cdot \mathbf w+b)\]
predicts whether the rank is $<2$ or $\ge 2$.

\subsection{Principal Component Analysis}\label{s:PCA}
Principal component analysis (PCA) is an unsupervised machine-learning strategy for dimensionality reduction.
In our case, we represent labelled elliptic curves as vectors in $\mathbb{R}^{1000}$, and PCA constructs a map $\mathbb{R}^{1000}\rightarrow\mathbb{R}^2$, the image of which groups curves according to their label.
The axes of this image are given by the principal components PC1 and PC2 of the dataset, from which the method takes its name.
A principal component is evaluated by $\sum_{n=1}^{1000} c_n a_{p_n}$ for $c_n \in \mathbb R$, and we call $c_n$ the {\em weight} of $a_{p_n}$ in the principal component.

\section{Experimental Results and Observations }\label{s:elliptic}
We now describe our new experimental results for elliptic curves defined over $\mathbb{Q}$. 
\subsection{Logistic regression for ranks 0, 1, and 2}\label{sec:ECrank}

Logistic regression is discussed in Section~\ref{s:logreg}.
In the previous paper \cite{HLOc}, it is demonstrated that logistic regression can distinguish elliptic curves of rank $0$ from those of rank $1$ with high accuracy based on their $a_p$-coefficients. 
With equation~\eqref{eq.parity}, that is, the parity conjecture, which is a theorem for curves of rank $r_E\in\{0,1\}$, in mind, one may wonder whether the classification is achieved through learning the parity of the root number $w(E)$. 
As mentioned in the Introduction, it is not completely clear how to extract the root number $w(E)$ from (a finite sequence of) the $a_p$-coefficients of an elliptic curve $E$.

To determine any possible role of the parity conjecture in the classification, we consider elliptic curves of rank $r_E \in\{0,1,2\}$, and perform logistic regressions for the datasets of elliptic curves with $r_E \in \{0,1\}, \{0,2 \}, \{1,2\}, \{0,1,2\}$, respectively. 
The curves in the datasets are all in the conductor range $[1, 1 \times 10^5]$, and a sample of $20,000$ curves are randomly chosen from each rank to make balanced datasets.

The results of our experiments are summarized in Table~\ref{t:rank}.
We see that all the accuracies are all over $0.96$, and sometimes over $0.99$. 
In particular, classification of rank 0 and 2 curves has accuracy $0.996$ and shows that it is not through recognizing the parity of $w(E)$. 
(Here we assume that the parity conjecture is true for rank $\le 2$.) The even higher accuracy $0.9998$ in the case of rank 1 and 2 tells us that rank 2 elliptic curves clearly distinguish themselves from those of rank 0 and 1 by way of the $a_p$-coefficients, and we will see it more clearly in the following subsections. 

\begin{table}[h!!!]
\begin{center}
{\footnotesize \begin{tabular}{|c|c|c|c|c|c|}
\hline
$N_E$ range & $r_E$ & |Data| = $\#\{E\}$&Precision&Confidence\\
\hline
$[1,1\times10^5]$   &$\{0, 1\}$  & 
                     $2.0\times10^4$ ($\times2$) & 0.961 & 0.92 \\
\hline
" &  \{0,2\}& " &  0.996& 0.99\\
\hline
" &  \{1,2\} &"&  0.999&0.99\\
\hline
" &  \{0, 1, 2\} & $2.0\times10^4$ ($\times3$) & 0.975 & 0.96 \\
\hline
\end{tabular}
}
\end{center}
\caption{{\sf
A table recording the results of logistic regression experiments for elliptic curve rank.
The first three rows after the header shows the precision and confidence (Matthews correlation coefficient) of a logistic regression binary classifier when asked to distinguish elliptic curves over $\mathbb{Q}$ by their rank. 
In the final row, we use a multinomial logistic regression classifier to distinguish all three ranks simultaneously. 
All of these experiments used a random sample of $2.0\times10^4$ curves for each rank, all with conductor $N_E$ in the range $[1, 1 \times 10^5]$.
}}
\label{t:rank}
\end{table}

\subsection{PCA for ranks 0, 1, and 2}\label{s:pca012}
\begin{figure}[h!!!]
\begin{center}
     \includegraphics[width=0.9\textwidth]{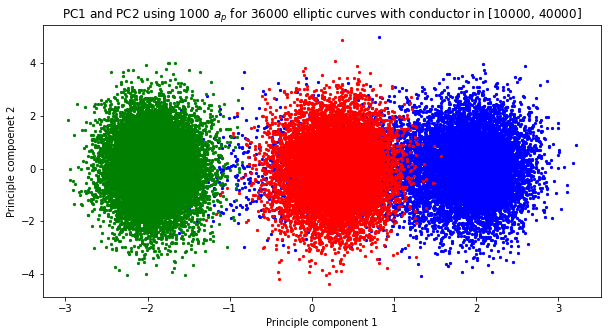}
\end{center}
\caption{\sf A plot of PC1 ($x$-axis) against PC2 ($y$-axis) for elliptic curves in the balanced dataset of $36,000$ randomly chosen elliptic curves with rank $r_E\in\{0,1,2\}$ and conductor $N_E\in[10000,40000]$. The blue (resp. red, green) points are the images of the vectors $v_L(E)$ corresponding to the elliptic curves in our dataset with rank $0$ (resp. $1$, $2$) under a map $\mathbb{R}^{1000}\rightarrow\mathbb{R}^2$ constructed using PCA. }
\label{pca-1}
\end{figure}

PCA is discussed in Section~\ref{s:PCA}.
We begin with PCA applied to a balanced dataset of $36,000$ randomly selected elliptic curves, each of which has rank $r_E\leq 2$ and conductor $N_E \in [1 \times 10^4, 4 \times 10^4]$. 
In particular, we plot these elliptic curves according to their PC1 and PC2 scores. 
The result is depicted in Figure~\ref{pca-1}, which shows a clear separation of all three ranks, with some overlap between rank 0 (blue) and rank 1 (red) curves.

\begin{figure}[h!!!]
\begin{center}
     \includegraphics[width=0.9\textwidth]{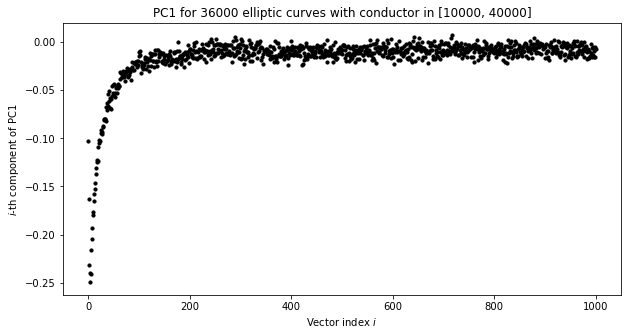}
\end{center}
\caption{\sf A plot of $n$ ($x$-axis) against the weight of $a_{p_n}$ ($y$-axis) in PC1 for the balanced dataset of $36,000$ randomly chosen elliptic curves with rank $r_E\in\{0,1,2\}$ and conductor $N_E\in[10000,40000]$ used in this Section.  }
\label{pca-2}
\end{figure}

We can better understand this separation by looking at the weights of $a_p$ in the principal components. Since the ranks are only separated along PC1, we will look only at the weights of this component.
It is clear from Figure \ref{pca-2} that the first hundred or so $a_p$ are the most important for this classification. Here we enumerate primes as $p_1=2, p_2=3, p_3=5, \dots$ and the $x$-axis represents the indices $1 \le n \le 1000$ for primes with $p_{1000}=7919$. 

\begin{figure}[h!!!]
\begin{center}
     \includegraphics[width=0.9\textwidth]{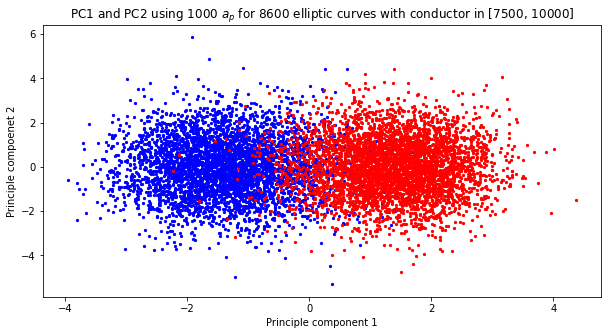}
\end{center}
\caption{\sf A plot of PC1 against PC2 for the balanced dataset of $8,600$ randomly chosen elliptic curves with rank $r_E\in\{0,1\}$ and conductor $N_E\in[7500,10000]$. The blue (resp. red) points correspond to curves with rank $0$ (resp. $1$).}
\label{pca-3}
\end{figure}

 We may further investigate the separation between rank 0 and rank 1 by removing the rank 2 data. 
In order to connect our observations to the oscillations in the average values of $a_p$ that have been presented in the introductio,n and are to be considered in the next subsection, we restrict the conductor interval to $[7500, 10000]$.
This reduces our dataset to just 4300 curves per rank.
The result of PCA is shown in Figure \ref{pca-3}.

\begin{figure}[h!!!]
\begin{center}
     \includegraphics[width=0.9\textwidth]{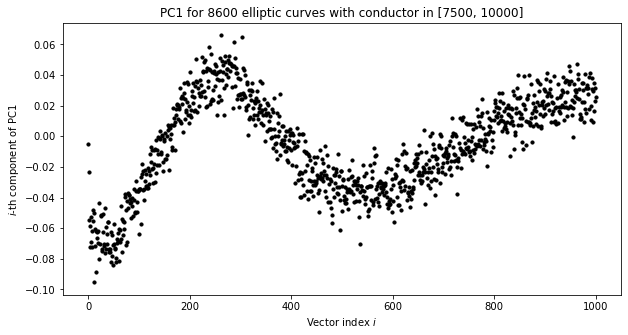}
\end{center}
\caption{\sf A plot of $n$ ($x$-axis) against the weight of $a_{p_n}$ in PC1 ($y$-axis) for the balanced dataset of $8,600$ randomly chosen elliptic curves with rank $r_E\in\{0,1\}$ and conductor $N_E\in[7500,10000]$.}
\label{pca-4}
\end{figure}

With this narrower conductor range, the separation of rank 0 and rank 1 according to PC1 is significantly better. We also see a fundamental difference in the PC1 from this dataset, as shown in Figure \ref{pca-4}. In particular, the weights of $a_p$ seem to follow a smooth, decaying oscillation. This indicates that there is an interesting structure in the datasets of $a_p$ which separates rank 0 and rank 1. In the following subsection, we will find such a structure in a statistical relationship between $p$ and $a_p$ for a fixed rank and conductor range. 

\subsection{Averages of $a_p$-coefficients}\label{sec:av}

In this subsection, we plot the averages of the $a_p$-coefficients to reveal surprising features that seem unknown in the literature. 

Let us begin by studying $f_0(n)$ and $f_1(n)$ for $N_E \in [7500, 10000]$, where $f_r(n)$ is defined as in equation \eqref{eq.frp}, which we repeat here for convenience: 
\begin{equation}\label{eq.frn-repeat}
f_r(n)=\frac{1}{\#\mathcal{E}_r[N_1,N_2]}\sum_{E\in\mathcal{E}_r[N_1,N_2]}a_{p_n}(E),
\end{equation}
where $p_n$ is the $n$th prime, and $N_1,N_2\in\mathbb{Z}_{>0}$ satisfy $N_1<N_2$, and $\mathcal{E}_r[N_1,N_2]$ is the set of (representatives of the isogeny classes of) elliptic curves over $\mathbb{Q}$ with rank $r$ and conductor in range $[N_1,N_2]$.
The expression in equation~\eqref{eq.frn-repeat} is simply the arithmetic mean of each $a_{p}$ over a set of elliptic curves with fixed rank and conductor range. 

\begin{figure}[h!!!]
\begin{center}
     \includegraphics[width=0.9\textwidth]{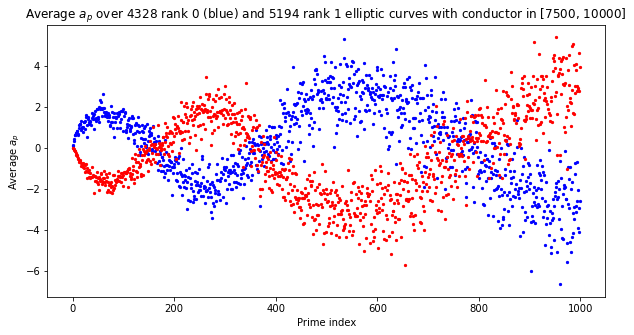}
\end{center}
\caption{\sf Plots of the functions $f_0(n)$ (blue) and $f_1(n)$ (red) for $1 \le n \le 1000$ and $[N_1,N_2]=[7500,10000]$.}  
\label{avg-1}
\end{figure}

In Figure~\ref{avg-1}, we present a plot for these functions. Observe that the values of $f_0(n)$ and $f_1(n)$ appear to follow an oscillation whose amplitude and period both grow with $n$. 
Moreover, these oscillations appear to mirror each other. 
Also note that the frequency of this oscillation matches that of the PC1 components in Figure \ref{pca-4} of the previous section in the sense that the oscillations attain $0$ approximately at the same $n$. 

We see a similar oscillation when looking at $f_2(n)$, although the pattern breaks for first several primes $p$. 
In Figure \ref{avg-2} below, we expand the conductor range to $[5000, 10000]$ in order to increase the number of rank 2 curves available. 
Even with this increase, we only have 1380 curves. The comparatively low number of curves likely contributes to the less concentrated distribution.
Note that including smaller conductor curves has slightly increased the frequency of oscillation. 
Indeed, we observe that as we look at elliptic curve sets of larger conductors, the frequency of oscillation becomes lower. We will further study the relationship between conductor and frequency in Section \ref{s:fitting}. 
\begin{figure}[h!!!]
\begin{center}
     \includegraphics[width=0.9\textwidth]{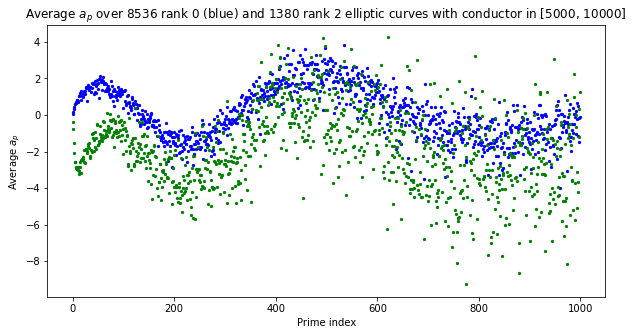}
\end{center}
\caption{\sf Plots of the functions $f_0(n)$ (blue) and $f_2(n)$ (green) for $1 \le n \le 1000$ and $[N_1,N_2]=[5000,10000]$.} 
\label{avg-2}
\end{figure}

\begin{figure}[h!!!]
\begin{center}
     \includegraphics[width=\textwidth]{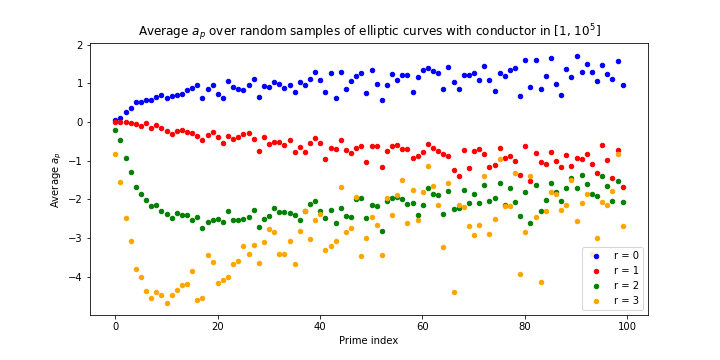}
\end{center}
\caption{\sf A plot of $f_r(n)$ for $r \in \{0,1,2,3\}$ and $N_E \in [1, 1 \times 10^5]$. The blue (resp. red, green, yellow) points correspond to curves of rank $0$ (resp. $1$, $2$, $3$). }
\label{avg-3}
\end{figure}

In Figure~\ref{avg-3}, we plot $f_r(n)$ for $r \in \{0,1,2,3\}$ and $N_E \in [1, 1 \times 10^5]$. 
While taking the average over such a large conductor range makes the oscillation much less apparent, it is worth noting that the average $a_{p_n}$ for $n \leq 20$ (i.e., $p \le 73$) distinguish all four ranks across this entire conductor range. Here, the yellow points correspond to average $a_p$ of rank $3$ curves. 
Also, note that we only have 531 rank 3 curves, whereas the other ranks are all plotted using a random sample of 20,000 curves.

\begin{center}
\end{center}

\subsection{Histograms of $a_p$ distributions}

\begin{figure}[h!!!]
\begin{center}
\includegraphics[width=.24\linewidth,valign=m]{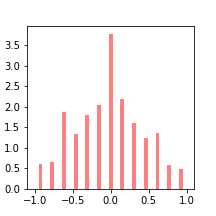} 
\includegraphics[width=.24\linewidth,valign=m]{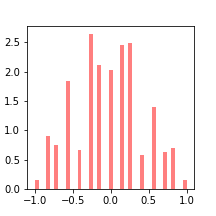}  \includegraphics[width=.24\linewidth,valign=m]{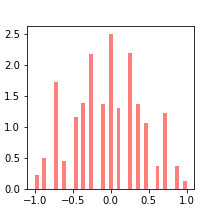}  \includegraphics[width=.24\linewidth,valign=m]{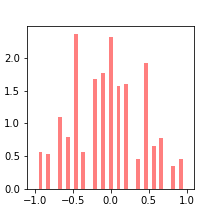} 
\end{center}
\caption{\sf (Left to right) Histograms of $\tilde{a}_{11}$, $\tilde{a}_{13}$, $\tilde{a}_{17}$, $\tilde{a}_{19}$
for curves of rank $0$.}
\label{t:rank0}
\end{figure}

\begin{figure}[h!!!]
\begin{center}
\includegraphics[width=.24\linewidth,valign=m]{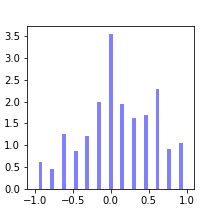} 
\includegraphics[width=.24\linewidth,valign=m]{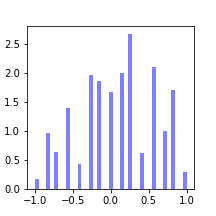}  \includegraphics[width=.24\linewidth,valign=m]{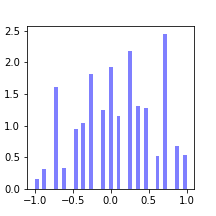}  \includegraphics[width=.24\linewidth,valign=m]{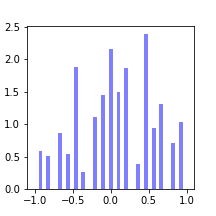}
\end{center}
\caption{\sf (Left to right) Histograms of $\tilde{a}_{11}$, $\tilde{a}_{13}$, $\tilde{a}_{17}$, $\tilde{a}_{19}$
for curves of rank $1$.}
\label{t:rank1}
\end{figure}

\begin{figure}[h!!!]
\begin{center}
\includegraphics[width=.24\linewidth,valign=m]{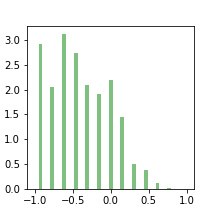} 
\includegraphics[width=.24\linewidth,valign=m]{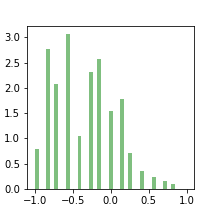}  \includegraphics[width=.24\linewidth,valign=m]{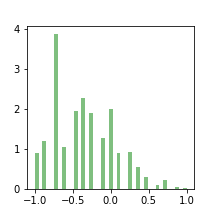}  \includegraphics[width=.24\linewidth,valign=m]{new_NormalizedR2a_17DensityPlot_7500,_10000_.png}
\end{center}
\caption{\sf (Left to right) Histograms of $\tilde{a}_{11}$, $\tilde{a}_{13}$, $\tilde{a}_{17}$, $\tilde{a}_{19}$
for curves of rank $2$.}
\label{t:rank2}
\end{figure}

In this section, we will look at the distribution of the normalized $a_p$ coefficient:
\begin{equation}\label{eq.normalized}
\tilde{a}_p=\frac{a_p}{2\sqrt{p}}
\end{equation}
for fixed $p$, and for elliptic curves with $N_E \in [7500, 10000]$ and rank $r_E\in\{0,1,2\}$. 
In equation~\eqref{eq.normalized} we normalize $a_p$ by the Hasse bound so that $\tilde{a}_p\in[-1, 1]$ for all $p$. 
In Figure \ref{t:rank0} (resp. \ref{t:rank1}, \ref{t:rank2}), we present the distributions of $\tilde{a}_p$ for curves of rank 0 (resp. 1, 2) and $p\in\{11,13,17,19\}$.

\begin{figure}[h!!!]
\begin{center}
\includegraphics[width=.24\linewidth,valign=m]{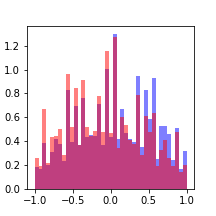} 
\includegraphics[width=.24\linewidth,valign=m]{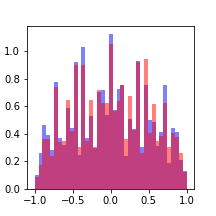}  \includegraphics[width=.24\linewidth,valign=m]{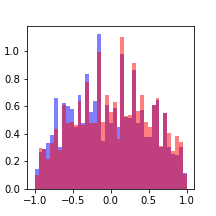} 
\includegraphics[width=.24\linewidth,valign=m]{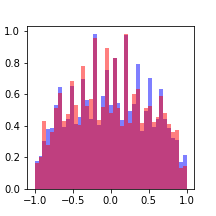}  
\end{center}
\caption{\sf (Left to right) Histograms of $\tilde{a}_{397}$, $\tilde{a}_{1151}$, $\tilde{a}_{1787}$, $\tilde{a}_{2731}$
for curves of rank $r_E\in\{0,1\}$.  In blue (resp. red), we have curves of rank $0$ (resp. 1), and, in purple, the distributions overlap.}
\label{t:skews}
\end{figure}

Whilst the rank 2 distributions are generally quite different, we note that the rank 0 and rank 1 distributions exhibit many similarities.
One exception is that rank 0 has a slight left skew and rank 1 has a slight right skew.
Based on the oscillation we see in Figure \ref{avg-1}, we expect that the rank 0 and rank 1 distributions will grow more skew until about $p = 397$, return to symmetric at about $p = 1151$, then become skewed in the opposite direction at around $p = 1787$, and finally return to symmetric at around $p = 2731$. 
In Figure~\ref{t:skews}, we present histograms of $\tilde{a}_p$, $p \in \{397, 1151, 1787, 2731\}$, for curves of rank 0 and rank 1  to see this phenomenon. 
Note that the purple shows where the distributions overlap.

We conclude this subsection by discussing the relationship these distributions have with the classification problem in Section~\ref{sec:ECrank}. 
For a fixed conductor range and a fixed $p$, these distributions are dependent on the rank of the curves from which they were generated. Moreover, each elliptic curve is associated to a specific infinite sequence of $(a_{p_n})_{n\geq1}$ which are essentially drawn at random from these distributions. 
By looking at sufficiently many of these draws, the classifier can predict with high accuracy whether they came from a sequence of rank 0, rank 1, or rank 2 distributions. 
The overwhelming overlap of the rank 0 and rank 1 distributions explain why this is the most difficult binary classification problem of the three, and why significantly more $a_p$ are required. 
In particular, as demonstrated in Section \ref{sec:heuristic}, using just $(a_{p_n})_{n=1}^{10}$, we can distinguish rank 2 from rank 0 or rank 1 with $\approx 0.97$ accuracy. 
In distinguishing rank 0 from rank 1, the accuracy is roughly 0.7 with 10 $a_p$-coefficients, depending on the sample.

\subsection{Curve fitting for the averages of $a_p$}\label{s:fitting}

Next, we turn to the question of curve fitting for the average $a_p$ plots presented in Section~\ref{sec:av}. 
Actually, we slightly modify equation~\eqref{eq.frp} and introduce the following function of primes $p \in\mathbb{Z}_{\geq1}$ (instead of the $n$th prime $p_n$):
\begin{equation}\label{eq.grp}
g_r(p)=\frac{1}{\#\mathcal{E}_r[N_1,N_2]}\sum_{E\in\mathcal{E}_r[N_1,N_2]}a_{p}(E),
\end{equation}
where $N_1 < N_2\in\mathbb{Z}_{>0}$, and $\mathcal{E}_r[N_1,N_2]$ is the set of (representatives of the isogeny classes of) elliptic curves over $\mathbb{Q}$ with rank $r$ and conductor in range $[N_1,N_2]$.
In particular, we want to find a curve which best approximates $g_0(p)$ and $g_1(p)$. 

\begin{figure}[h!!!]
\begin{center}
     \includegraphics[width=\textwidth]{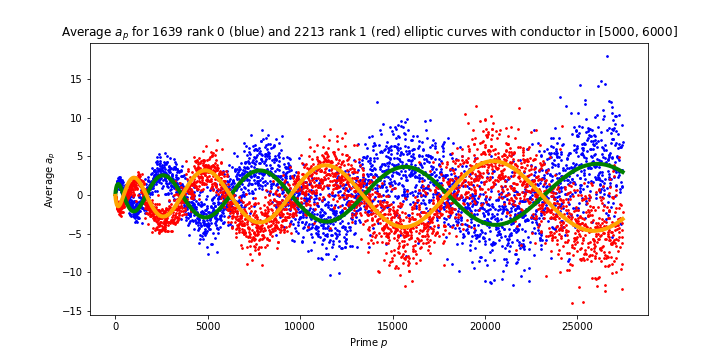}
\end{center}
\caption{\sf Plot of $g_0(p)$ (res. $g_1(p)$) in blue (resp. red) for elliptic curves with conductor $N_E\in[5000,6000]$. The solid curve is the corresponding curve of best fit $0.5398x^{0.1980}\sin({0.1255x^{0.5272}})$ (resp. $-0.4875x^{0.2215}\sin({0.1239x^{0.5288}})$), which has mean squared error $7.4768$ (resp. $5.4605$).}
\label{f-cf1}
\end{figure}

\begin{figure}[h!!!]
\begin{center}
     \includegraphics[width=\textwidth]{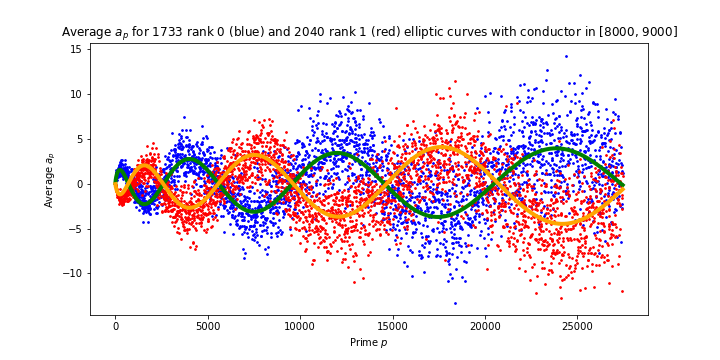}
\end{center}
\caption{\sf Plot of $g_0(p)$ (res. $g_1(p)$) in blue (resp. red) for elliptic curves with conductor $N_E\in[8000,9000]$. The solid curve is the corresponding curve of best fit $0.5243x^{0.2004}\sin({0.0948x^{0.5331}})$ (resp. $-0.2581x^{0.2828}\sin({0.0994x^{0.5277}})$), which has mean squared error $7.3044$ (resp. $6.4967$).}
\label{f-cf2}
\end{figure}

\begin{figure}[h!!!]
\begin{center}
     \includegraphics[width=\textwidth]{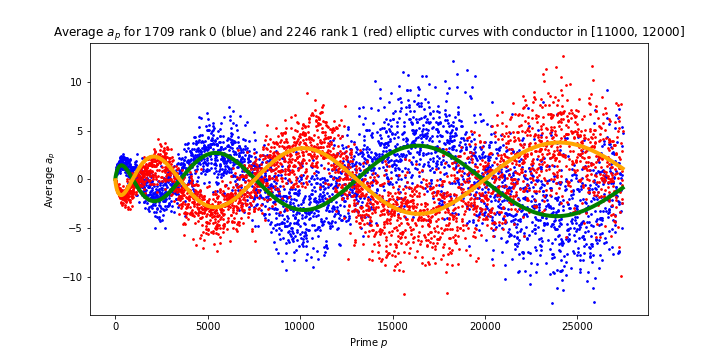}
\end{center}
\caption{\sf Plot of $g_0(p)$ (res. $g_1(p)$) in blue (resp. red) for elliptic curves with conductor $N_E\in[11000,12000]$. The solid curve is the corresponding curve of best fit $0.4273x^{0.2160}\sin({0.0835x^{0.5291}})$ (resp. $-0.5400x^{0.1934}\sin({0.0871x^{0.5246}})$), which has mean squared error $7.3127$ (resp. $5.8253$).}
\label{f-cf3}
\end{figure}

\begin{figure}[h!!!]
\begin{center}
     \includegraphics[width=\textwidth]{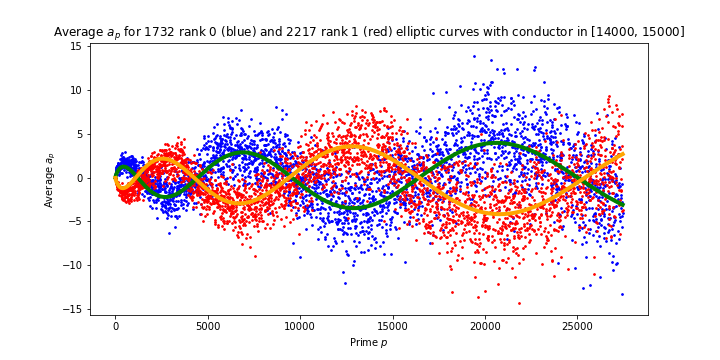}
\end{center}
\caption{\sf 
Plot of $g_0(p)$ (res. $g_1(p)$) in blue (resp. red) for elliptic curves with conductor $N_E\in[14000,15000]$. The solid curve is the corresponding curve of best fit $0.2273x^{0.2884}\sin({0.0727x^{0.5308}})$ (resp. $-0.2013x^{0.3048}\sin({0.0863x^{0.5131}})$), which has mean squared error $7.1645$ (resp. $6.09544$).}
\label{f-cf4}
\end{figure}

Plotting points $(p,g_r(p))$ yields an oscillation with increasing amplitude and period, and so trigonometric polynomial of a linear argument seems inadequate. 
Motivated by this observation, we will look at curves of the form 
\begin{equation}\label{eq.bfs}
y=Ax^{\alpha}\sin\left(Bx^{\beta}\right),
\end{equation}
where the parameters $A, \alpha, B, \beta$ are tuned to minimize the mean squared error. 
In Figure~\ref{f-cf1} (resp.~\ref{f-cf2},~\ref{f-cf3}, \ref{f-cf4}), we plot $g_0(p)$ and $g_1(p)$ along with the curve of best fit for conductor range $[5000,6000]$ (resp. $[8000,9000]$, $[11000, 12000]$, $[14000,15000]$).

Note that the aforementioned relationship between conductor and frequency can be clearly observed both in these pictures and in the parameter $B$. 
However, it is difficult to make this relationship more rigorous when $\beta$, which also effects frequency, is not held constant. 
It is also worth noting that $\beta \approx 0.5$ in all the fits we have tried.
In Table~\ref{t:bfs}, we record numerical values for $(A,\alpha,B,\beta)$ and mean squared errors (MSEs) to two decimal places for several conductor ranges.

\begin{table}[h!!!]
\begin{center}
{\footnotesize 
\begin{tabular}{| *{5}{c|}}
\hline
& \multicolumn{2}{c|}{$r_E=0$} &\multicolumn{2}{c|}{$r_E=1$}\\
\hline
$N_E$ range & $(A,\alpha,B,\beta)$ &MSE& $(A,\alpha,B,\beta)$ & MSE \\
\hline
$[5000,6000]$&$(0.54,0.20,0.13,0.53)$ & $7.48$&$(-0.47,0.22,0.12,0.53)$ & $5.46$\\
$[6000,7000]$&$(0.74,0.16,0.11,0.53)$& $7.80$&$(-0.50,0.21,0.11,0.53)$ & $5.51$\\
$[7000,8000]$&$(0.55,0.19,0.11,0.52)$ & $7.46$&$(-0.39,0.24,0.11,0.52)$ & $6.30$\\
$[8000,9000]$& $(0.52,0.20,0.09,0.53)$ & $7.30$& $(-0.26, 0.28, 0.10,0.53)$ & $6.50$\\
$[9000,10000]$ & $(0.36,0.24,0.08,0.54)$ & $6.96$& $(-0.30,0.27, 0.09,0.53)$ & $6.29$ \\
$[10000,11000]$ & $(0.32,0.25,0.08,0.54)$  & $7.83$& $(-0.46,0.21,0.08,0.54)$ & $5.61$\\
$[11000,12000]$ & $(0.43,0.22,0.08,0.53)$  & $7.31$& $(-0.54, 0.19,0.09,0.52)$ & $5.83$\\
$[12000,13000]$ & $(0.42,0.22,0.09,0.52)$  & $7.49$& $(-0.39,0.23,0.09,0.51)$ & $5.75$\\
$[13000,14000]$ & $(0.55,0.19,0.08,0.52)$  & $7.03$& $(-0.31,0.23,0.09,0.51)$ & $6.12$\\
$[14000,15000]$ & $(0.23,0.29,0.07,0.53)$  & $7.16$& $(-0.20,0.30,0.09,0.51)$ & $6.10$\\
\hline
\end{tabular}
}
\end{center}
\caption{{\sf
Table recording best fit parameters for the functions $g_0(p)$ and $g_1(p)$ and the mean squared errors (MSEs) for elliptic curves in the specified conductor range. Numerical values are rounded to two decimal places.
}}
\label{t:bfs}
\end{table}

\begin{remark}
The function $g_r(p)$ depends not only on $p$, but also the conductor range $[N_1,N_2]$ and the rank $r$. 
Looking at Table~\ref{t:bfs}, it seems that the error does not vary much with the conductor range, but does vary with the rank (and rank 1 errors seem to be smaller).
\end{remark}

\begin{remark}
On first glance, it might appear that equation \eqref{eq.bfs} suggests the average value of $a_p$ grows like $O(p^{\alpha})$, and Table~\ref{t:bfs} suggests that $\alpha$ is very roughly $1/5$. 
The Ramanujan conjecture, which is a theorem in this case, implies that $a_p=O(p^{1/2})$.
In the absence of bounds for the error, we make no precise claims about the average growth of $a_p$.
\end{remark}


\section{Heuristic Classification}\label{sec:heuristic}

In this section we will present a heuristic function for binary classification of elliptic curves of rank $\leq1$ and of rank $\ge 2$ for a fixed conductor range. 
The function will be {\em heuristic} in the sense that it is a simple function which approximates the classification.  More precisely, the value of the function represents the probability of an elliptic having rank $\ge 2$ as a result of machine learning. Applying a threshold of $0.5$, this leads to a binary classification that achieves high accuracy.

We offer some motivation for the heuristic function developed here as follows.
The vast majority of elliptic curves have rank equal to $0$ or $1$.
Empirically, on the LMFDB, of the $16494$ curves with conductor between $7500$ and $10000$, there are $15538$ with rank $0$ or $1$.
In an asymptotic sense, it is conjectured that 50\% of all elliptic curves over $\mathbb{Q}$ have rank $0$ and 50\% have rank $1$ (for further discussion of the origins of this conjecture, see \cite[Introduction]{BS1}).
It might therefore be natural to separate the task of rank classification into first distinguishing ranks $\geq2$ from those $\leq1$.
Since our heuristic function is very simple, it can be readily used to get the probability for a given curve to have rank $\ge 2$. 
An alternative approach to classifying higher rank curves was developed in \cite{KV}.

From Figure \ref{avg-3}, it is clear that the first several $a_p$-coefficients (on average) distinguish the case of rank $\leq1$ from that of rank $\ge 2$. 
This observation suggests that we use only the first 10 or so $a_p$-coefficients to perform logistic regression for fixed conductor ranges. 
Indeed, using just $(a_{p_1},a_{p_2},\dots,a_{p_{10}})=(a_2,a_3,\dots,a_{29})$, we were able to distinguish curves of rank $<2$ from curves of rank $\geq 2$ in the conductor range $[10000, 20000]$ with accuracy $\approx 0.95$. 
This was done on a balanced dataset of 3400 curves of each class. 
By reducing the conductor range to $[1, 10000]$ and to a dataset of 1970 curves of each class, we obtained accuracy $\approx 0.97$. 

For $\mathbf a:=(a_2, a_3, a_5, \dots ,  a_{29})$, a heuristic function for conductor range [10000,20000] can be defined by
\[ r(\mathbf a) = \frac 1 {1 + e^{- (\mathbf w \cdot \mathbf a + b)}} , \]
where $(\mathbf w; b)=(w_1, w_2, \dots , w_{10}; b)$ is determined by logistic regression and has entries 
\begin{multline*}
 \left(-1.1198144 , -1.12733444 , -0.98921727, -0.87923555, -0.57809252,-0.51279302, \right.\\
 \left.-0.32884407, -0.3539072,  -0.24136925, -0.19393439;  -5.62771169\right).
\end{multline*}
Here $b$ corresponds to a bias and we use $0.5$ as a threshold. 

\begin{example}
Consider the elliptic curve given by $$y^2+xy+y=x^3-150508x+13027931.$$ This curve has rank $1$, conductor $15015$ and $\mathbf a=(1, 1, 1, 1, -1, 1, 2, -4, -4, 6)$. 
Then we obtain $r(\mathbf a)= 0.00011 <0.5$. 
On the other hand, if we consider the elliptic curve given by $$y^2=x^3+x^2-436x-336,$$ it has rank $2$, conductor $15080$ and $\mathbf a = (0, -2, -1, -4, 0, 1, 0, -8, 0, -1)$. 
For this curve, we get $r(\mathbf a) = 0.97456 >0.5$. 
In this way, we can distinguish two cases heuristically with high accuracy.
\end{example}

For the conductor range $[1,10000]$, the vector $(\mathbf w;b)$ can be taken to be
\begin{multline*}
 \left(-1.41299148, -1.77879752, -1.38817256, -1.03428287, -0.71286324, -0.59119957,\right.\\
 \left. -0.40613106, -0.39675042, -0.2878296, -0.22388697; -8.77332846\right).
\end{multline*}
 In \cite{AP}, the fourth author offered a more systematic investigation into  the relationship between a conductor range and a number of $a_p$-coefficients required to attain high accuracy in heuristic classification. 


\section{Conclusions and Outlook}\label{s:outlook}

In practical terms, the experimental results presented in this article further demonstrate the utility of data-scientific approaches in arithmetic classification problems.
Of course, these experiments can be generalised in several directions, for example, replacing elliptic curves by rational modular forms of weight $>2$ or arithmetic curves of genus $>1$; replacing the base field $\mathbb{Q}$ by number fields of larger degree; or replacing the rank by other invariants of interest such as the Tate--Shafarevich group order.
The successful implementation of basic machine-learning strategies is a continuation of the theme developed in \cite{HLOa}, \cite{HLOb}, \cite{HLOc}.

On the other hand, one might seek to generalise the methodology, for example, incorporating more unsupervised machine-learning techniques, or possibly applying methods of reinforcement learning.
Whilst the application of modern approaches into a classical subject such as number theory is interesting, perhaps more exciting in this article is the appearance of seemingly new mathematical structures within the data.

In particular, we highlight the unexpected and striking behaviour of function $f_r(n)$ introduced in equation~\eqref{eq.frp}.
In Section~\ref{s:pca012}, it was noted that PCA can shed some light through weights on the oscillations observed.
That said, several immediate questions remain. 
For example, one might seek to quantify the error of the approximations developed in Sections~\ref{s:fitting} in terms of $p$, and subsequently explore new implications for the variation of $a_p$ in families of elliptic curves. 
Furthermore, whilst the Sato--Tate conjecture asserts that the average value of $a_p$ cannot grow monotonically with $p$, it is completely mysterious that the value should oscillate in such a notable manner.
These questions, and several others, may be formulated in mathematical terms, but the answers may involve some interplay with machine-learning techniques.


{\small 
Yang-Hui He {\sf hey@maths.ox.ac.uk} \\
London Institute for Mathematical Sciences, Royal Institution, London W1S 4BS, UK \\
Department of Mathematics, City, University of London, EC1V 0HB, UK

Kyu-Hwan Lee {\sf khlee@math.uconn.edu} \\
Department of Mathematics, University of Connecticut, Storrs, CT 06269-1009, USA \\
Korea Institute for Advanced Study, Seoul 02455, Republic of Korea

Thomas Oliver {\sf T.Oliver@westminster.ac.uk} \\
University of Westminster, London, UK

Alexey Pozdnyakov {\sf alexey.pozdnyakov@uconn.edu} \\
Department of Mathematics, University of Connecticut, Storrs, CT 06269-1009, USA
}

\end{document}